\documentclass[a4paper, 12pt]{article}
\usepackage{enumerate, theorem}
\usepackage{amsmath, amsfonts, amssymb}
\usepackage[height=22.5cm, width=15cm]{geometry}
\usepackage[all]{xy}
\usepackage{mathrsfs, yfonts}

\DeclareMathOperator{\id}{id}

\DeclareMathOperator{\C}{\mathbb{C}}

\newcommand{\parag}[1]{\paragraph{\sc{#1.}} }

\newtheorem{thm}{Theorem}[subsection]
\newtheorem{defn}[thm]{Definition}
\newtheorem{cor}[thm]{Corollary}
\newtheorem{prop}[thm]{Proposition}
\newtheorem{lemma}[thm]{Lemma}

\begin{document}

\title{New singularity invariants : the sheaf $\beta_{X}^{\bullet}$. }

 \author{Daniel Barlet\footnote{Institut Elie Cartan, G\'eom\`{e}trie,\newline
Universit\'e de Lorraine, CNRS UMR 7502   and  Institut Universitaire de France.}.}

\maketitle

\tableofcontents

\parag{Abstract}The graded coherent sheaf $\alpha_{X}^{\bullet}$ constructed in [B.18] for any reduced pure dimensional complex space $X$ is stable by exterior product but not by the de Rham differential. We construct here a new graded coherent sheaf $\beta_{X}^{\bullet}$ containing $\alpha_{X}^{\bullet}$ and stable both by exterior product and  by the de Rham differential. We show that it  has again the ``pull-back property''  for  holomorphic maps $f : X \to Y$ between irreducible complex spaces such that $f(X)$ is not contained in the singular set of $Y$. Moreover, this graded coherent sheaf $\beta_{X}^{\bullet}$ comes with a natural coherent exhaustive filtration  and this filtration is also compatible with the  pull-back by such holomorphic maps. These sheaves define  new invariants on singular complex spaces. We show on some simple examples that these invariants are new.

\parag{AMS classification 2020} 32 S 05 - 32 S 10 - 32 S 20.

\section{Erratum for ``the sheaf $\alpha_{X}^{\bullet}$''}

The aim of this section is to correct several mistakes in [B.18]. The main mistake is in theorem 4.1.1 which is wrong in the very general setting in which it is stated.\\

So we begin with a much more modest version of the ``pull-back`` theorem'' for these sheaves which has a rather simple proof.

\begin{thm}\label{pull-back}
Let $f : X  \to Y$ be a holomorphic map between irreducible complex spaces and assume that $f(X)$ is not contained in the singular set $S(Y)$ of $Y$. Then there exists a natural ``pull-back map''
$$\hat{f}^{*} : f^{*}(\alpha_{Y}^{\bullet}) \to \alpha_{X}^{\bullet} $$
which extends the usual pull-back of the sheaf $f^{*} : f^{*}(\Omega_{Y}^{\bullet}\big/torsion) \to \Omega_{X}^{\bullet}\big/torsion$.\\
For any holomorphic maps $f : X \to Y$ and $g : Y \to Z$ between irreducible complex spaces such that $f(X) \not\subset S(Y)$ and  $g(f(X)) \not\subset S(Z)$  we have
$$\hat{f}^{*}(g^{*}(\sigma)) = \widehat{f\circ g}^{*}(\sigma) \quad \forall \sigma \in \alpha_{Z}^{\bullet} .$$
\end{thm}

\parag{Proof} The problem is local. Let $\sigma$ be a  section of the sheaf $\alpha_{Y}^{\bullet}$ on an open set $V$ in $Y$. Let $V'$ be the set of regular points in $V$ and let $U''$ the set of regular points in the open set $U' := f^{-1}(V')$. This is a Zariski dense open set in $U := f^{-1}(V)$ and, as $\sigma$ is a holomorphic form on $V'$, $f^{*}(\sigma)$ is a well defined holomorphic form on $U''$ which is Zariski open and dense in $U$. Take a point $x$ in $U$; by definition (see proposition 2.2.4 in \cite{[B.18]}) there exists a open neighborhood $W$ of $y := f(x)$ in $V$ and a monic polynomial
$$P(z) = z^{k} + \sum_{h=1}^{k} S_{h}.z^{k-h} $$
such that $S_{h}$ is a section on $W$ of the symmetric algebra of degree $h$, $S_{h}(\Omega_{Y}^{\bullet}\big/torsion)$, of the sheaf $\Omega_{Y}^{\bullet}\big/torsion$, which satisfies  $P(\sigma) = 0$ in $\Gamma(W, S_{k}(\Omega_{Y}^{\bullet}\big/torsion))$. Then the pull-back $f^{*}(P)$ of $P$ by $f$ is well defined on $f^{-1}(W)$ and is a monic polynomial whose coefficients are sections on $f^{-1}(W)$ of the symmetric algebra of $\Omega_{X}^{\bullet}\big/torsion$. On the open set $U''\cap f^{-1}(W)$ the holomorphic form $f^{*}(\sigma)$ is a root of $f^{*}(P)$ and so the meromorphic\footnote{Remember that $\sigma$ is a meromorphic form on $V$ with poles in $S(Y)\cap V$.} form $f^{*}(\sigma)$ on $U\cap f^{-1}(W)$ is integrally dependent on the sheaf $\Omega_{X}^{\bullet}\big/torsion$. So it defines a unique section on $U$ of the sheaf $\alpha_{X}^{\bullet}$.\\
The second assertion of the theorem is a simple consequence of the fact that the sheaf $\alpha_{X}^{\bullet}$ has no torsion.$\hfill \blacksquare$\\

The second mistake (which is consequence of the previous one) is that in definition 5.1.5 of \cite{[B.18]} it is necessary to ask that the $p-$dimensional irreducible analytic subset  $Y$ is not contained in the singular set of $X$ to define the integral on $Y$ of a form of the type $\rho.\alpha\wedge \bar \beta$, where $\alpha, \beta$ are sections of the sheaf $\alpha_{X}^{p}$ in $X$.\\

To be clear we give here the correct statements for definition 5.1.5, lemma 5.1.6 and for the theorem 5.1.7. The proof  given in \cite{[B.18]} of this theorem is correct but makes sense only assuming that the pull-back for the sheaf $\alpha^{\bullet}$ are defined. This is consequence of the following hypotheses which allow to apply the corrected version of the theorem 4.1.1 above.

\begin{defn}\label{integral 0}
Let $X$ be an irreducible  complex space and let $Y \subset X$ be a closed irreducible $p-$dimensional analytic subset in $X$; assume that $Y$ is not contained in the singular set  $S(X)$ of $X$. We shall note $j : Y \to X$ the inclusion map. Let $\rho$ be a continuous function with compact support in $X$ and let $\alpha, \beta$ be sections on $X$ of the sheaf $\alpha_{X}^{p}$. We define the number $\int_{Y} \ \rho.\alpha\wedge \bar \beta$ as the integral
$$ \int_{Y} \ j^{*}(\rho).\hat{j}^{*}(\alpha)\wedge \overline{\hat{j}^{*}(\beta)} $$
which is well defined by corollary 5.1.2 of \cite{[B.18]} using the inclusion $\alpha_{Y}^{p} \subset L_{Y}^{p}$.
\end{defn}

\bigskip

\begin{lemma}\label{integral 2}
Let $f : X \to Y $ be a holomorphic map between irreducible complex spaces Let $Z$ be a closed $p-$dimensional irreducible analytic subset in $X$ such that $Z$ is not contained in the singular set $S(X)$ of $X$,  the restriction of $f$ to $Z$ is proper and $f(Z)$ is not contained in the singular set of $Y$.
Let $\alpha, \beta$ be sections on $Y$ of the sheaf $\alpha_{Y}^{p}$ and let $\rho$ be a continuous function with compact support in $Y$. Then we have the equality
$$ \int_{Z} \ f^{*}(\rho).\hat{f}^{*}(\alpha)\wedge \overline{\hat{f}^{*}(\beta)} = \int_{f_{*}(Z)} \ \rho.\alpha \wedge \bar \beta .$$
If $f(Z)$ is contained in $S(Y)$ the singular set of $Y$ and has dimension at most $p-1$ we have $ \int_{Z} \ f^{*}(\rho).\hat{f}^{*}(\alpha)\wedge \overline{\hat{f}^{*}(\beta)} = 0$.
\end{lemma}

\parag{Proof} The first assertion is an  easy consequence of the same result when $\alpha, \beta$ are holomorphic forms, by considering a modification of $Z$ where it is the case, using for instance, a desingularization of $Z$ (see \cite{[H.64]} normalizing the sheaf $\Omega_{Z}^{p}$ (see \cite{[B.18]}).\\
 When $f(Z) \subset S(Y)$ and $f_{*}(Z) = 0$ the restriction of $f$ to $Z$ has generic rank at most $p-1$, so the pull-back of any holomorphic $p-$form on $Y$ to $Z$ is $0$. Then the monic polynomial giving an integral dependence relation of $\alpha$ (or of $\beta)$) reduces to $z^{k} = 0$ on $f(Z)$ and so $\alpha$ (and $\beta$) vanishes on $Z$. $\hfill \blacksquare$

\bigskip

\begin{thm}\label{integral 3}
Let $X$ be an irreducible complex space and $(Y_{t})_{t \in T}$ be an analytic family of $p-$cycles in $X$ parametrized by a reduced complex space $T$. Assume that for $t$ in a dense open subset $T'$ in $T$ no component of the cycle $Y_{t}$ is contained in $S(X)$, the singular set of $X$. Let $\rho$ be a continuous function with support  in the compact set  $K$  in $X$ and let $\alpha, \beta$ be two sections of the sheaf $\alpha_{X}^{p}$. Define the function
$$ \varphi : T' \to \C \quad {\rm by} \quad \varphi(t) := \int_{Y_{t}} \ \rho.\alpha\wedge \bar \beta .$$
Then $\varphi$ is continuous on $T'$ and locally bounded near each point in $T$. More precisely, for any continuous hermitian metric $h$ on $X$ and  any compact set $L$ in $T$, there exists a constant $C > 0$ (depending on $K, \alpha, \beta, h, L$ but not of the choice of  $\rho$ with support in $K$) such that for each $t \in T'$ we have:
$$ \vert \varphi(t)\vert \leq C.\int_{Y_{t}} \ \vert\rho\vert.h^{\wedge p} \leq C.\vert\vert \rho\vert\vert.\int_{Y_{t} \cap K} \ h^{\wedge p}.$$
\end{thm}

Remark that the subset of points $t \in T$ where the cycle $Y_{t}$ has no irreducible component meeting a given compact set and contained in $S(X)$ is a closed analytic subset in $T$ by a general result on analytic families of cycles (see for instance the proposition IV 3.5.7 in \cite{[B-M 1]} in the case of compact cycles). So, assuming that $T$ is irreducible, if there exists a point $t$ such that $Y_{t}$ has no irreducible component meeting $K$ and contained in $S(X)$, there exists  a Zariski open and dense subset $T'$ of $T$ which satisfies the hypothesis in the previous theorem.\\

The last mistake is the lemma 6.2.2 which is wrong for $k \geq 4$. The correct computation of $\alpha_{S_{k}}^{2}$ is given in the  paragraph 3.3.

\section{Definition of $\beta_{X}^{\bullet}$ and the pull-back property}

The next sections of this paper are  complements to \cite{[B.18]}. So the notations are the same than in {\it loc. cit.} In particular, for a reduced pure dimensional complex space $X$ we note $L_{X}^{\bullet}$ the graded sheaf of meromorphic forms on $X$ which are holomorphic on any desingularization of $X$ and $\omega_{X}^{\bullet}$ the sheaf of $\bar\partial-$closed currents on $X$ of type $(\bullet, 0)$ modulo its torsion sub-sheaf. Recall that the graded sheaf $\alpha_{X}^{\bullet}$ constructed in {\it loc. cit.} is a graded coherent  sub-sheaf of $L_{X}^{\bullet}$ which is again  a graded coherent sub-sheaf of $\omega_{X}^{\bullet}$. Of course, all these sheaves contain the sheaf $\Omega_{X}^{\bullet}\big/torsion$  and  coincide with it on the non singular part  of $X$.\\

\subsection{Construction of the sheaf $\beta_{M}^{\bullet}$}

Let $X$ be a pure dimensional  reduced complex space and note $\alpha_{X}^{\bullet}$ the graded sheaf on $X$ introduced in \cite{[B.18]}.

\begin{lemma}\label{exterior product}
The sheaf $\alpha_X$ is stable by exterior product.
\end{lemma}

\parag{Proof} Recall that, by definition, the sheaf $\alpha_ X^{\bullet}$ is a sub-sheaf of the sheaf $\omega_ X^{\bullet}$ and a section $\sigma$ on the open set $U \subset X$ of the sheaf $\omega_ X^{\bullet}$ is a section on $U$ of $\alpha_ X^{\bullet}$ if it may be written locally on $U$   as $\sigma = \sum_{j \in J} \ \rho_{j}.\omega_{j}$ where $\omega_j$ are holomorphic forms on $U$ and $\rho_j$ are $\mathscr{C}^{\infty}$ functions on the complement of the singular set $S$ in $U$ which are  bounded near $S$. Is it clear that the exterior product of two such sections on $U$ of $\alpha^{\bullet}_ X$ can be written in the same way locally on $U$ and then define a current on $U$ which is $\bar\partial-$closed on $U \setminus S$. So to conclude the lemma, it is enough to prove that this current admits a $\bar \partial-$closed extension to $U$.  In fact, as  the sheaf $\alpha^{\bullet}_ X$ is a sub-sheaf of the sheaf $L^{\bullet}_ X$ obtained by the direct image of the sheaf $\Omega^{\bullet}_{\tilde{ X}}$ where $\tau : \tilde{ X} \to  X$ is any desingularisation of $X$  and as this sheaf  $L_{X}^{\bullet}$ is stable by exterior product, the conclusion follows from the inclusion $L_{X}^{\bullet} \subset \omega_{X}^{\bullet}$.$\hfill \blacksquare$\\

Remark that the sheaf $\alpha^{\bullet}_ X$ is a graded $\Omega^{\bullet}_ X-$module but is not stable in general by the de Rham differential. For instance in $ X := \{(x, y, z) \in \C^3 \ / \  x.y = z^2\}$ the differential form $ dx\wedge dy/z = - d\big(z.dx/x - z.dy/y)$ is not in $\alpha^2_ X$ but the form $z.dx/x - z.dy/y$ is a section of $\alpha^1_ X$ (see \cite{[B.18]} or paragraph 3.3).

\parag{A construction} Define $\alpha^{\bullet}_ X[0] := \alpha^{\bullet}_ X$ and for any integer $p \geq 0$ and any integer $q \geq 0$ define 
\begin{equation}
 \alpha^{q}_ X[p+1] := \sum_{r=0}^q \Big( \alpha_ X^r[p] \wedge \alpha^{q-r}_ X[p]\Big) +   \sum_{r=0}^{q-1} \Big(\alpha_ X^r[p] \wedge d\big( \alpha^{q-r-1}_ X[p] \big)\Big)\subset L^q_ X
 \end{equation}
Recall that the sheaf $L_{X}^{\bullet}$ is stable by exterior products and by the de Rham differential.\\
\begin{prop}\label{construction}
Then we have the following properties:
\begin{enumerate}
\item For each integer  $p$ the sheaf $\alpha_{X}^{\bullet}[p]$ is stable by exterior product with $\Omega_{X}^{\bullet}\big/torsion$. Moreover  for each integers $p, q$ we have $\alpha_{X}^{0}.\alpha_{X}^{q}[p] \subset \alpha_{M}^{q}[p] $.
\item For each integers $p, q$ the sheaf $\alpha^q_ X[p]$ is $\mathcal{O}_ X-$coherent sub-sheaf of $L^q_ X$.
\item For each integers $p, q$ the sub-sheaf $\alpha^q_ X[p]$ is contained in $\alpha^q_ X[p+1]$.
\item For each integers  $p, q$ and $ q'$ we have $\alpha^{q}_ X[p] \wedge \alpha^{q'}_ X[p] \subset \alpha^{q+q'}_ X[p+1]$.
\item For each integers $p, q$ and $ r$ we have $\alpha^r_ X[p]\wedge d\big(\alpha^q_ X[p]\big) \subset \alpha^{q+r+1}_ X[p+1]$.\\
 In particular  $d\big(\alpha^q_ X[p]\big) \subset \alpha^{q+1}_ X[p+1]$.
\end{enumerate}
\end{prop}

\parag{Proof} The  property 1. is an obvious consequence of the definition of these sheaves by an induction on $p$.\\
As $L_{ X}^{q}$ is a coherent sheaf on $ X$, to prove 2. it is enough to prove that $\alpha_{ X}^{q}[p+1]$ is a finite type $\mathcal{O}_{ X}-$module. We shall prove this by an induction on $p \geq 0$.\\
So assumed that for each  $q$ the coherence of the sheaf $\alpha_{X}^{q}[p]$. Then we want to prove that $\alpha_{X}^{q}[p+1]$ is finitely generated. Let $g_{j,r}$ be a generator of the sheaf $\alpha_{X}^{r}[p]$. Then we shall show that the elements $ g_{i,r}\wedge g_{j, q-r}$ and $ g_{i, r}\wedge dg_{j, q-r-1}$ for all choices of $i$ and $j$, generates $\alpha_{X}^{q}[p+1]$. The only point which is not obvious is the fact that for any sections $u \in \alpha_{X}^{r}[p]$ and $v \in \alpha_{X}^{q-r-1}[p]$ the wedge product $u \wedge dv$ is in the sheaf generates by our ``candidates'' generators. But then write
\begin{align*}
& u = \sum_{i} a_{i}.g_{i, r} \quad {\rm and} \\
& v = \sum_{j} b_{j}.g_{j, q-r-1}
\end{align*}
where $a_{i}$ and $b_{j}$ are holomorphic functions. Then
$$dv = \sum_{j} db_{j}\wedge g_{j, q-r-1} + \sum_{j} b_{j}.dg_{j, q-r-1} .$$
So in the wedge products $u\wedge dv$ the terms are linear combinations of our candidates generators excepted those like $a_{i}. g_{i,r}\wedge db_{j}\wedge g_{j, q-r-1}$. This  point is solved by the condition 1. which is already proved. 
The points  3. 4. and 5. are obvious. This complete the proof of our induction.$\hfill \blacksquare$\\

Now remark that the sequence of coherent sub-sheaves $\alpha_{X}^{\bullet}[p]$ of the coherent sheaf $L_{X}^{\bullet}$ is increasing. So it is locally stationary on $X$.

\begin{defn}\label{beta}
Define the coherent sub-sheaf $\beta_{ X}^{\bullet} $ as the union  of the increasing sequence of coherent sub-sheaves $\alpha_{ X}^{\bullet}[p], p \geq 0$ of the coherent sheaf $L_{ X}^{\bullet}$.
\end{defn}

\begin{cor}\label{obvious 1}
The graded sub-sheaf $\beta_{ X}^{\bullet}$ of  the graded coherent  differential  sheaf $L_{ X}^{\bullet}$ is coherent, stable by exterior product and by the de Rham differential.
\end{cor}

\parag{Proof} The assertion is local, so we may assume that $\beta_{X}^{\bullet} = \alpha_{X}^{\bullet}[p] \quad \forall  p \geq p_{0}$. Then the proposition is consequence of the properties $3.$ and $4.$ above.$\hfill \blacksquare$\\

\begin{thm}\label{not obvious 2}
For any holomorphic map $f : X \to Y$ between irreducible complex spaces such that $f(X)$ is not contained in the singular set $S(Y)$ of $Y$,  there exists an unique  pull-back
$$ \hat{f}^{*}:  f^{*}(\beta_{Y}^{\bullet}) \to \beta_{X}^{\bullet} $$
which is compatible with the pull-back of the $\alpha^{\bullet}-$sheaves (see section 1) and which is graded of degree $0$ and compatible with the exterior product and the  de Rham differential.
For any holomorphic maps $f : X \to Y$ and $g : Y \to Z$ between irreducible complex spaces, such that $f(X) \not\subset S(Y)$ and $g(f(X)) \not\subset S(Z)$ we have
$$ \widehat{f\circ g}^{*}(\sigma) = \hat{f}^{*}(\hat{g}^{*}(\sigma)) \qquad \forall \sigma \in \beta_{Z}^{\bullet} .$$
Moreover, for each integer $p \geq 0$ the pull-back $\hat{f}^{*}$  induces a pull-back
 $$\hat{f}^{*}[p] : f^{*}(\alpha_{Y}^{\bullet}[p]) \to \alpha_{X}^{\bullet}[p]$$
  and in the previous situation $ \widehat{f\circ g}^{*}[p] = \hat{g}^{*}[p]\circ \hat{f}^{*}[p]$ for each $p \geq 0$.
\end{thm}

So we  shall construct in fact a (graded) naturally filtered sheaf $(\beta_{X}^{\bullet}, (\alpha_{X}^{\bullet}[p])_{p \in \mathbb{N}})$ such that the  pull-back constructed in the previous theorem is compatible with these filtrations and with the composition of suitable holomorphic maps.\\

In the following  we  make the convention that $\alpha_{X}^{\bullet}[-1] := \alpha_{X}^{\bullet}[0] := \alpha_{X}^{\bullet} $.\\

\parag{Proof} Assume that $X$ is an irreducible complex space and that $f(X)$ is not contained in the singular locus $S(Y)$ of $Y$.  Assume also that for some integer $p \geq 0$ we have constructed for any $q \leq p$ a pull-back morphism
$$ \hat{f}_{q}^{*} : f^{*}(\alpha_{Y}^{\bullet}[q]) \to \alpha_{X}^{\bullet}[q] $$
with the following properties
\begin{enumerate}[$1_{p}$]
\item It induces the usual pull-back of the sheaves of holomorphic forms when it is restricted to the smooth parts\footnote{Our hypothesis implies that there exists a dense Zariski open set $X''$ in $X \setminus S(X)$ such that the restriction of $f$ to $X''$ takes values in $Y \setminus S(Y)$.} of $X$ and $Y$. Note that this implies that the restriction of $\hat{f}_{p}^{*}$ to $ f^{*}(\alpha_{Y}^{\bullet}[q])$ is equal to $\hat{f}_{q}^{*}$ because, by definition the sections of the sheaves under consideration are determined by their restrictions to an open dense subset.
\item For $s, t$ in $\alpha_{Y}^{\bullet}[p-1]$ we have $\hat{f}_{p}^{*}(s \wedge t) = \hat{f}_{p-1}^{*}(s) \wedge \hat{f}_{p-1}^{*}(t)$.
\item For any $u$ in  $\alpha_{Y}^{\bullet}[p-1]$ such that $du$ is in $\alpha_{Y}^{\bullet+1}[p]$\footnote{For $p \geq 1$ $u \in \alpha_{X}^{\bullet}[p-1]$ implies $du \in \alpha_{X}^{\bullet+1}[p]$ is automatic; but not for $p = 0$ with our convention.} we have $d(\hat{f}_{p-1}^{*}(u)) = \hat{f}_{p}^{*}(du)$.
\end{enumerate}
Then we want to construct $\hat{f}_{p+1}^{*} : f^{*}(\alpha_{N}^{\bullet}[p+1]) \to \alpha_{M}^{\bullet}[p+1] $ satisfying again the  properties above for $p+1$.\\

It is clear that  that our inductive hypothesis given by the conditions $1_{p}, 2_{p}$ and $3_{p}$  is true for $p = 0$ (but $2_{0}$ is obtain by looking at points in $X''$ and using the absence of torsion). \\
Now we shall show that if it is satisfied for some $p \geq 0$ then it is also satisfied for $p + 1$.

\parag{Construction of $\hat{f}_{p+1}^{*}$} Let $\xi$ be a section in $\alpha_{Y}^{\bullet}[p+1]$. We may write
$$ \xi = \sum_{j= 0}^{J} \beta_{j}\wedge \gamma_{j} + \sum_{j=0}^{J} u_{j}\wedge dv_{j} $$
where $\beta_{j}, \gamma_{j}, u_{j}, v_{j}$ are sections of the sheaf $\alpha_{Y}^{\bullet}[p]$. It is clear that our conditions $1_{p+1}, 2_{p+1}, 3_{p+1}$ implies that we must put
$$\hat{f}_{p+1}^{*}(\xi) = \sum_{j=0}^{J} \hat{f}_{p}^{*}(\beta_{j})\wedge \hat{f}_{p}^{*}(\gamma_{j}) + \sum_{j=0}^{J} \hat{f}_{p}^{*}(u_{j})\wedge d( \hat{f}_{p}^{*}(v_{j})) .$$
Now the main point is to prove that if we change the choice of writing $\xi$ in such a way, the value of $\hat{f}_{p+1}^{*}(\xi)$ stays the same. In other words, we have to prove that if $\xi = 0$ is written as above then we find $\hat{f}_{p+1}^{*}(\xi) = 0$.\\
To prove this is quite simple because it is enough to look on $X''$. On this open dense subset we have simply taken the usual pull-back of the holomorphic form $\xi$ restricted to the smooth part of $Y$ by the holomorphic map $f' : X'' \to Y\setminus S(Y)$ induced by $f$. As this pull-back commutes with exterior product and de Rham differential, its result is independent on the way with have written $\xi$ above. This implies our claim because the sheaf $\alpha_{X}^{\bullet}[p+1]$ has no torsion.\\
To verify the properties $1_{p+1}, 2_{p+1}$ and $3_{p+1}$ is then obvious because it is enough to check them on $X''$.\\
 This completes the proof of the existence of pull-back morphisms $\hat{f}^{*}[p]$ for each $p \geq 0$ and then for the sheaves $\beta^{\bullet}$. And it also gives the compatibility of these pull-back with the exterior product and the de Rham differential.\\
 The only point which we have to precise to complete the proof of the theorem \ref{not obvious 2} is the ``functorial'' aspect of these pull-back. But this is again an easy consequence of  the non existence of torsion for the sheaves we consider. $\hfill \blacksquare$\\

\subsection{Filtration}

\begin{prop}\label{filt.}
Let $X$ be an irreducible complex space. Then for each $q \geq 0$ we have $\beta_{X}^{q} = \alpha_{X}^{q}[q]$. If $X$ is normal, for $q \geq 1$ we have $\beta_{X}^{q} = \alpha_{X}^{q}[q-1]$.
\end{prop}

\parag{Proof} First remark that, by definition $\alpha_{X}^{0}$ is the sheaf of locally bounded meromorphic functions on $X$ (so it is  equal to $\mathcal{O}_{X}$ if and only if $X$ is normal), and that $\beta_{M}^{0} = \alpha_{M}^{0}$ by definition.\\
Fix an integer $q_{0} \geq 1$ and assume that for any integer $ q < q_{0}$ we have $\alpha_{X}^{q}[p-1] = \beta_{X}^{q}$ for some integer $p \geq 1$. This means that $\alpha_{X}^{q}[p-1] = \alpha_{X}^{q}[p] $ for these $q$. By definition we have
$$ \alpha_{X}^{q_{0}}[p+1] = \sum_{h=0}^{q_{0}} \alpha_{X}^{h}[p]\wedge \alpha_{X}^{q_{0}-h}[p] + \sum_{h=0}^{q_{0}-1} d(\alpha_{X}^{h}[p])\wedge \alpha_{X}^{q_{0}-h-1}[p] .$$
But our assumption allows to replace $p$ by $p-1$ in the right hand-side of the equality above, so we find that $\alpha_{X}^{q_{0}}[p+1] = \alpha_{X}^{q_{0}}[p] $.\\
Now remark that,  for each $q \geq 0$, $\alpha_{X}^{q}$ is stable by multiplication by elements in $\alpha_{X}^{0}$, so we have
$$ \alpha_{X}^{1}[1] = \alpha_{X}^{1}[0] + \sum_{i, j=1}^{I} \mathcal{O}_{X}.g_{j}.dg_{i} $$
where $g_{1}, \dots, g_{I}$ generate the coherent $\mathcal{O}_{X}-$module $\alpha_{X}^{0}$. This implies that $\alpha_{X}^{1}[1]$ is stable by multiplication by $\alpha_{X}^{0}$ and this implies the equality
$$ \alpha_{X}^{1}[2] = \alpha_{X}^{1}[1] + \sum_{i, j=1}^{I} \mathcal{O}_{X}.g_{j}.dg_{i} = \alpha_{X}^{1}[1].$$
So we have $\alpha_{X}^{1}[1] = \beta_{X}^{1}$. This allows to begin our induction on $q_{0}$ for $q_{0} = 1$ with $p = 2$.\\
Then by induction on $q_{0} \geq 1$ we conclude that for each $q \geq 1$ we have $\beta_{X}^{q} = \alpha_{X}^{q}[q]$.\\
 In the case where $X$ is normal, we may take $I = \{1\}$ and $g_{1} = 1$ and this shows that $\alpha_{X}^{1}[0] = \beta_{X}^{1}$ and the induction gives now, if we begin with $q_{0} = 1$ and $p = 1$, the equality   $\alpha_{X}^{q}[q-1] = \beta_{X}^{q}$ for each $q \geq 1$.$\hfill \blacksquare$\\

\parag{Remark} This shows that for a normal complex space we always have the equality $\beta_{X}^{1} = \alpha_{X}^{1}$, so the sheaf $\beta_{X}^{\bullet}$ is ``new'' only in degrees at least equal to $2$ when $X$ is normal.

\subsection{Integration and Stokes' formula}

We begin by defining the integration on $p-$cycles.

\begin{lemma}\label{integration}
Let $X$ be an irreducible  complex space and $u, v$ be sections of the sheaf $\beta_{X}^{p}$. Let $Y$ be an irreducible $p-$cycle in $X$ which is not contained in $S(X)$ and let $\rho$ be a continuous function with compact support on $ X$. Then the improper integral
$$ \int_{Y} \ \rho.\hat{j}(u)\wedge \overline{\hat{j}(v)} $$
is absolutely convergent, where $j : Y \to X$ is the inclusion map.
\end{lemma}

\parag{Proof} This is an easy consequence of the fact that, using the pull-back theorem, on a suitable modification of $Y$ the $p-$forms $u$ and $v$ becomes holomorphic.$\hfill \blacksquare$\\

Remark that for any holomorphic function $f$ on $Y$ which does not vanish on a non empty open set in $Y$ this integral is the limit as $\varepsilon > 0$ does to zero of the integral
$$ \int_{Y \cap \{ \vert f \vert > \varepsilon \}} \ \rho.\hat{j}(u)\wedge \overline{\hat{j}(v)}.$$

\begin{defn}\label{integration bis}
In the situation of the previous lemma, the number
 $$ \int_{X} \ \rho.\hat{j}(u)\wedge \overline{\hat{j}(v)} $$
 will be called the  integral of $ \rho.\hat{j}(u)\wedge \overline{\hat{j}(v)}$ on $Y$. This definition extends by linearity to general $p-$cycles in $X$ which have no irreducible component in $S(X)$.
\end{defn}

In the sequel we shall write simply $\int_{Y} \ \rho.u\wedge \bar v$ this integral, omitting the pull-back by $j$. But it is necessary to keep in mind that this abuse of notation is acceptable because of the compatibility of the pull-back maps for the sheaf $\beta$ with the composition of suitable maps.

\begin{lemma}\label{direct image}
Let $X$ be an irreducible complex space and $f : X \to Y$ an holomorphic map. Let $Z$ be an irreducible $p-$cycle in $X$ such that its direct image $f_{*}(Z)$ is defined (as a $p-$cycle\footnote{Remember that this means that the restriction of $f$ to $Z$ is proper.}) in $Z$. Assume that $Z$ is not contained in $S(X)$ and that $f(Z)$ is not contained in $S(Y)$. Let $u, v$ be sections of the sheaf $\beta_{Y}^{p}$ and $\rho$ a continuous function on $Y$ with compact support. Then we have
$$ \int_{f_{*}(Z)} \ \rho.u\wedge \bar v = \int_{Z} \ f^{*}(\rho).\hat{f}(u) \wedge \overline{\hat{f}(v)}.$$
\end{lemma}

\parag{Proof} Thanks to the lemma \ref{integration} and the remark following it, this reduces to the usual change of variable in the case where $f(Z)$ is not contained in $S(Y)$.$\hfill \blacksquare$\\

Again the functoriality of the pull-back implies here the fact that we may either take the pull-back by $f$ and restrict to $ Z$ or directly take the pull-back by the restriction of $f$ to $Z$.

\begin{prop}\label{integration 2}
Let $(Y_{t})_{t \in T}$ be an analytic family of $p-$cycles in $X$ parametrized by a  reduced complex space $T$ and let $\rho :  X \to \C$ be a continuous function with compact support in $X$. Assume that for each $t \in T$  the cycle $Y_{t}$ has no irreducible component contained in $S(X)$. Let $u, v$ be sections of the sheaf $\beta_{X}^{p}$. Then the function $\varphi : T \to \C$ defined by
$$ \varphi(t) = \int_{Y_{t}} \ \rho.u\wedge \bar v  \quad \forall t \in T $$
is continuous on a dense open set $T'$ of $T$ and is locally bounded near each point in $T$.
\end{prop}

\parag{Proof} Consider a modification $\tau : \tilde{X} \to X$ of $X$ with center in $S(X)$ such that $u$ and $v$ becomes holomorphic on it. Let $\nu : \tilde{T} \to T$ be the  normalization of  $T$.

\parag{Claim} There exists a modification $\sigma : \Theta \to \tilde{T}$ with $\Theta$ a normal complex space and  an analytic family $(\tilde{Y}_{\theta})_{\theta \in\Theta}$ of $p-$cycles  in $\tilde{X}$ such that for $\theta$ generic, $\tilde{Y}_{\theta}$ is the strict transform of $Y_{\nu(\sigma(\theta))}$ by $\tau$. Moreover, for each $\theta \in \Theta$ we will have $\tau_{*}(\tilde{Y}_{\theta}) = Y_{\nu(\sigma(\theta))}$. \\

\parag{Proof of the claim} As $\tilde{T}$ is normal, we may decompose in an open neighborhood of the support of $\rho$ the family $(Y_{\tilde{t}})_{\tilde{t} \in \tilde{T}}$ as a finite sum of analytic families of $p-$cycles  parametrized by $\tilde{T}$ such that their generic cycles are irreducible. Then, by additivity of the integral, it is enough to treat the case of such a family. This means that , without no loss of generality, we may assume that the family $Y_{\tilde{t}})_{\tilde{t} \in \tilde{T}}$ has an irreducible generic cycle.\\
So let $G \subset \tilde{T} \times X$ the graph of the family ; it is irreducible. Let $\Gamma$ the strict transform of $G$ by the modification $\id_{\tilde{T}}\times \tau : \tilde{T}\times \tilde{X} \to \tilde{T} \times X$.\\
Now we may find a modification $\sigma : \Theta \to \tilde{T}$ such that the strict transform $\tilde{\pi} : \tilde{\Gamma} \to \Theta$ by the modification $\sigma $ of the projection $\pi : \Gamma \to \tilde{T}$ becomes equidimensionnal. So, as $\Theta$ is normal, the fibres of $\tilde{\pi}$ give an analytic family of $p-$cycles in $\tilde{X}$ parametrized by $\Theta$, such that for $\theta$ generic, $\tilde{Y}_{\theta}$ is the strict transform of $Y_{\nu(\sigma(\theta))}$. Moreover, for each $\theta \in \Theta$ we will have $\tau_{*}(\tilde{Y}_{\theta}) = Y_{\nu(\sigma(\theta))}$ because this is true for $\theta$ generic and both are analytic families of cycles in $X$, thanks to the direct image theorem (see \cite{[B-M 1]} theorem IV 3.5.1). This proves the claim.\\
Then the result follows using the previous lemma and the continuity of the integration of a  continuous form with compact support on a continuous family of cycles (see \cite{[B-M 1]} proposition IV 2.3.1).$\hfill \blacksquare$\\

\parag{Remark} The only point where we use the fact that $u, v$ are sections of the sheaf $\beta_{X}^{p}$ is when we define the integral using the definition \ref{integration} and the second part of the  theorem \ref{not obvious 2} which gives that the pull-back is compatible with composition of suitable holomorphic maps.  In the rest of the proof, we only use the fact that $u, v$ are sections of $L_{X}^{p}$ to know that on a suitable modification of $X$ they become holomorphic forms.\\

\begin{lemma}\label{Stokes}
Let $X$ be a reduced complex space and let $u$ and $v$ be respectively sections of the sheaves $\beta_{X}^{p-1}$ and $\beta_{X}^{p}$. Let $\rho$ be a $\mathscr{C}^{1}-$function on $X$ with compact support. Then for any $p-$cycle $Z$ in $X$  which has no irreducible component in $S(X)$ we have
$$ \int_{X} \ d(\rho.u\wedge \bar v) = 0 .$$
\end{lemma}

\parag{Proof} It is enough to prove this formula when $Z$ is  irreducible. As $Z$ is not contained in the singular locus of $X$  the result follows from the fact that on a suitable modification of $Z$ the forms $u$ and $v$ becomes holomorphic so we may apply the standard Stokes's theorem.$\hfill \blacksquare$\\

Note that the commutation of the pull-back maps with the de Rham differential is crucial here.\\

\section{Examples}

\subsection{The sheaf $\alpha$ for a product}

\begin{prop}\label{product}
Let  $X$ and $Y$  be irreducible complex spaces. Then if $p_{1} : X \times Y \to X$ and $p_{2} : X\times Y \to Y$ are the projections, we have for each $p \geq 0$ a natural isomorphism
$$\theta : \oplus_{q = 0}^{p} \  p_{1}^{*}(\alpha_{X}^{q})\otimes_{\mathcal{O}_{X\times Y}} p_{2}^{*}(\alpha_{Y}^{p-q}) \to \alpha_{X\times Y}^{p} $$
given by exterior product.
\end{prop}

\parag{Proof} This is an easy exercice using two desingularizations $\sigma : \tilde{X} \to X$ and $\tau : \tilde{Y} \to Y$ which are normalizing respectively for the sheaves $\Omega_{X}^{\bullet}\big/torsion$ and $\Omega_{Y}^{\bullet}\big/torsion$, as the product map $\sigma\times \tau : \tilde{X}\times  \tilde{Y}  \to X\times Y$ is a desingularization of $X\times Y$ which normalizes the sheaf  $\Omega_{X\times Y}^{\bullet} $.
$\hfill \blacksquare$\\

\parag{Remark} It is easy to extend this proposition to the sheaves $\alpha_{X\times Y}^{\bullet}[p]$ for each $p \geq 0$ and so to the sheaf $\beta_{X\times Y}^{\bullet}$.\\

The following trivial corollary will be used in an example below.

\begin{cor}\label{utile}
Let $X$ be an irreduciblel complex spaces. Consider on  $X \times D$ a $L^{p+1}$ form   $\omega\wedge f(z).dz$  where $D$ is a disc in $\C$ with coordinate $z$, $f : D \to \C$ a holomorphic function on $D$ which is not identically zero,  and where $\omega$ is a $L^{p}-$form on $X$. Then $\omega$ is a section of $\alpha_{X}^{p}$ if and only if $\omega\wedge f(z).dz$ is a section of $\alpha_{X\times D}^{p+1}$.\end{cor}

\subsection{The curve $X := \{ x^{3} = y^{5}\} \subset \C^{2}$}

\begin{lemma}\label{curve}
On the curve  $X :=  \{ x^{3} = y^{5}\} \subset \C^{2}$ we have
\begin{align*}
& \alpha_{X}^{0} = L_{X}^{0} = \mathcal{O}_{X} \oplus \C.y^{2}/x \oplus \C.y^{4}/x^{2} \oplus \C.y^{3}/x  \oplus \C.y^{4}/x  \\
& \omega^0_{X}  = L_{X}^{0} +  \mathcal{O}_{X}.y/x^{2} \\
& \alpha_{X}^{1} = \Omega_{X}^{1}\big/torsion + \mathcal{O}_{X}.y^{2}.dy/x \\
& \beta_{X}^{1} = L_{X}^{1} =  \alpha_{X}^{1}[1] =  \mathcal{O}_{X}.y^{2}.dx/x^{2} \\
& \omega_{X}^{1} =  \Omega_{X}^{1}\big/torsion + \mathcal{O}_{X}.dy/x^{2}
\end{align*}
\end{lemma}

The proof is left as an exercise to the reader. For the computation of the sheaf  $\omega_{X}^{0}$  the reader may  use the fact that 
$$(y/x^{2}).(3x^{2} - 5y^{4}.dy)\big/(x^{3}- y^{5}) = (3y.dx - 5x.dy)\big/(x^{3}- y^{5})$$
 in $H_{X}^{1}(\Omega_{\C^{2}}^{1})$ and the characterization given in \cite{[B.78]} of the sheaf $\omega_{X}^{\bullet} $ in terms of the fundamental class of $X$.\\

 \subsection{The surfaces $S_{k}$} 
  Consider the surfaces  $S_{k}:= \{(x, y, z) \in \C^{3}\ / \ x.y = z^{k}\}$ for $k$ an integer at least equal to $2$.\\
 In the following lemma, we determine the sheaves $\alpha_{S_{k}}^{\bullet} $ and $\beta_{S_{k}}^{\bullet}$.  We also correct the  lemma 6.2.2 of [B.18] which is wrong for $k \geq 4$.
 
 \begin{lemma}\label{Sk}
Let $m := [k/2]$ be the integral part of $k/2$. Then we have 
 \begin{align*}
 & \beta_{S_{k}}^{1} = \alpha_{S_{k}}^{1} = \Omega_{S_{k}}^{1}\big/torsion + \mathcal{O}_{S_{k}}.x.dy\big/z^{m}  \\
 & \alpha_{S_{k}}^{2} =  \Omega_{S_{k}}^{2}\big/torsion  +  \mathcal{O}_{S_{k}}.\frac{dx\wedge dy}{z^{m-1}} \quad  \\
  & \beta_{S_{k}}^{2} = \alpha_{S_{k}}^{2}[1] =   \Omega_{S_{k}}^{2}\big/torsion + \mathcal{O}_{S_{k}}.\frac{dx\wedge dy}{z^{m}}.
  \end{align*}
 \end{lemma}
 
 \parag{Proof} The first assertion is consequence of the equality $\alpha_{M}^{1} = \beta_{M}^{1}$ for any normal complex space which is proved in proposition \ref{filt.}. The computation of $\alpha_{S_{k}}^{1}  $ is an obvious consequence of lemma 6.2.3  in \cite{[B.18]}. Note that the equalities
  $$x.dy/z^{m} + y.dx/z^{m}  = k.z^{k-m}.dz \quad {\rm and} \quad  (x.dy/z^{m}).( y.dx/z^{m}) = z^{k-2m}.(dx).(dy)$$
   gives the integral dependance relation of $x.dy/z^{m}$ on $S^{\bullet}(\Omega_{S_{k}}^{1}\big/torsion)$. \\
 Let now prove the second assertion.\\
 Remark first that we have on $S_{k}$ the relations
 $$ x.dx\wedge dy =  k.z^{k-1}.dx\wedge dz  \quad  y.dx\wedge dy = k.z^{k-1}.dz\wedge dy  $$
 and using the equality $x.y = z^{k}$ this implies
 $$ dx\wedge dy = k.y.dx\wedge dz/z = -k.x.dy\wedge dz/z .$$
 Dividing by $z^{m-1}$ this gives
 $$ \big(\frac{dx \wedge dy}{z^{m-1}}\big)^{2} = -k^{2}.z^{k-2m}.(dx\wedge dz).(dy\wedge dz) \quad {\rm in} \ S^{2}(\Omega_{S_{k}}^{2}\big/torsion).$$
 This prove that $dx \wedge dy/z^{m-1}$ is a section of the sheaf $\alpha_{S_{k}}^{2}$.\\
 
 We want to prove now that the meromorphic form
$$\frac{ dx\wedge dy}{z^{m}} = k.y.\frac{dx\wedge dz}{z^{m+1}} = - k.x.\frac{dy\wedge dz}{z^{m+1}} $$
 which corresponds to $k^{2}.(a.b)^{k-m}.da\wedge db$ via the quotient map
 $$ q_{k} : \C^{2} \to S_{k} \quad (a, b) \mapsto (x = a^{k}, y= b^{k}, z = a.b) $$
 is not in $\alpha_{S_{k}}^{2} = \Omega_{S_{k}}^{2}\big/torsion$.\\
 As the fiber $F_{0}$  of the sheaf  $ F := q_{k}^{-1}(\Omega_{S_{k}}^{2}\big/torsion)$ at $0$ is the $A := \C\{a^{k}, b^{k}, a.b\}-$sub-module of $A.da\wedge db$ generated by $a^{k}.da\wedge db, b^{k}.da\wedge db, (a.b)^{k-1}.da\wedge db$, we have to show that $(a.b)^{k-m-1}.da\wedge db$ is not integral on $F_{0}$. This an easy consequence of the fact that for $q < k/2$ there is no positive constant $C$ such that for $a > 0$ and $b > 0$ small enough we have the inequality $ (a.b)^{q } \leq C.(a^{k} + b^{k})$.\\

 To prove the last assertion remark first that the form $d(x.dy\big/z^{m})$ is in $\alpha_{S_{k}}^{2}[1] = \beta_{S_{k}}^{2}$ (this last equality is also proved in proposition \ref{filt.}). But we have on $S_{k}$, using the equality $[(k-1)/2] + [k/2] = k-1$
 \begin{align*}
 & y.dx + x.dy = k.z^{k-1}.dz \quad \quad  {\rm so} \\
 & x.dy\wedge dx = k.z^{k-1}.dz\wedge dx \quad {\rm and \  then} \\
 & \frac{dy\wedge dx}{z^{m}} = k.y.\frac{dz\wedge dx}{z^{m+1}}
 \end{align*}
 This gives $d(x.dy\big/z^{m}) = (1 - m/k).dx\wedge dy\big/z^{m}$.\\
  So the inclusion of  $ \Omega_{S_{k}}^{2}\big/torsion + \mathcal{O}_{S_{k}}.dx\wedge dy\big/z^{m}$ in $\beta_{S_{k}}^{2}$ is proved. The equality $\alpha_{S_{k}}^{2}[1] = \beta_{S_{k}}^{2}$ easily implies the equality in the previous inclusion, as we have the inclusion  $\alpha_{S_{k}}^{1}\wedge \alpha_{S_{k}}^{1} \subset \Omega_{S_{k}}^{2}\big/torsion$..$\hfill \blacksquare$\\

Note that the sheaf $L_{S_{k}}^{2}$ is equal to $\Omega_{S_{k}}^{2}\big/torsion + \mathcal{O}_{S_{k}}.dx\wedge dy\big/z^{k-1}$. So for $k \geq 4$ we have strict inclusions between $ \Omega_{S_{k}}^{2}\big/torsion ,  \alpha_{S_{k}}^{2}, \beta_{S_{k}}^{2}$ and $L_{S_{k}}^{2} = \omega_{S_{k}}^{2}$.

\subsection{$M_{k} := \{x.y = u^{k}.v \}$}

Let $m := [k/2]$ be the integral part of the integer $k \geq 1$.

\begin{lemma}\label{ex.1, 1}
The meromorphic $1-$form  $\omega_{m} := x.dy/u^{m}$ belongs to $\alpha_{M_{k}}^{1}$ but for $k \geq 2$ the differential  $d\omega_{m}$ is not in $\alpha_{M_{k}}^{2}$.
\end{lemma}

\parag{Proof} We have
\begin{align*} 
& x.dy/u^{m} + y.dx/u^{m} = d(xy)/u^{m} = d(u^{k}v)/u^{m} = k.u^{k-1-m}v.du + u^{k-m}.dv \quad {\rm and} \\
&  \big(x.dy/u^{m}\big).\big(y.dx/u^{m}\big) = x.y.(dx).(dy)/u^{2m} = u^{k-2m}.v.(dx).(dy) 
\end{align*}
so $\omega_{m}$ satisfies the following integral dependance relation on $\Omega_{M}^{1}\big/torsion$
\begin{equation}
\omega_{m}^{2} - (k.u^{k-m-1}.v.du + u^{k-m}.dv).\omega_{m} + u^{k-2m}.v.(dx).(dy) = 0 
\end{equation}
Now we have 
$$ d\omega_{m} = \frac{dx\wedge dy}{u^{m}} - m.\frac{x.du\wedge dy}{u^{m+1}} .$$
But now we restrict this $2-$form to the surface $S_{k} := \{v = 1\}$ which cuts the singular set of $M_{k}$ only at the point $x = y = u = 0, v=1$, and we find, as we have on this surface 
 $x.dy + y.dx = k.u^{k-1}.du$ which  implies $y.dx\wedge dy = k.u^{k-1}.du\wedge dy$ and then
   $u.dx\wedge dy = k.x.du\wedge dy $. So 
$$  (d\omega_{m})_{\vert \{ v= 1\}} = (1 - m/k). dx\wedge dy/u^{m} $$
which is not in $\alpha_{S_{k}}^{2}$ for $k \geq 2$ (see   lemma \ref{Sk}).$\hfill \blacksquare$\\


\begin{lemma}\label{ex.1, 2}
The $2-$form $w := \omega_{m}\wedge dv$ belongs to $\alpha_{M_{k}}^{2}$ but $dw$ is not in $\alpha_{M_{k}}^{3}$ for $k \geq 2$.
\end{lemma}

\parag{Proof} The first assertion is obvious as $\alpha_{M_{k}}^{\bullet}$ is stable by wedge products and contains $\Omega_{M_{k}}^{\bullet}\big/torsion$.\\
To prove the second assertion consider the following holomorphic map

$$ \pi : S_{k}\times \C \to M_{k}, \quad ((x, y, u)), v) \mapsto (x.v, y, u, v) .$$

Then $\pi^{*}(dw) = dx\wedge dy\wedge dv/u^{m} - m.x.du\wedge dy\wedge dv/u^{m+1}$. Using the corollary \ref{utile} of the  proposition \ref{product} and the fact that we have on $S_{k}\times \C$
$$\pi^{*}(dw) = v.dv \wedge \big((k-m).x.du\wedge dy/u^{m+1}\big) $$
we conclude that $\pi^{*}(dw) $ is not a section of $\alpha_{S_{k}\times \C}^{3}$,  concluding the proof.$\hfill \blacksquare$\\

\begin{cor}\label{Mk}
For $k \geq 4 $ we have on $M_{k}$
\begin{align*}
& \Omega_{M_{k}}^{1}\subset \alpha_{M_{k}}^{1} = \beta_{M_{k}}^{1} \subset L_{M_{k}}^{1} \\
& \Omega_{M_{k}}^{2}\subset \alpha_{M_{k}}^{2} \subset \beta_{M_{k}}^{2} \subset L_{M_{k}}^{2} \\
& \Omega_{M_{k}}^{3}\subset \alpha_{M_{k}}^{3} \subset  \beta_{M_{k}}^{3} \subset L_{M_{k}}^{3}
\end{align*}
where all inclusions are strict.
\end{cor}

We leave to the reader the easy proof using the previous computations.$\hfill \blacksquare$\\

\subsection{Fermat surfaces}

We shall look now to the surfaces $F_{n} := \{ (a, b, z) \in \C^{3} \ / \  a^{n} - b^{n} = z^{n}\}$ for $n \geq 3$. As these surfaces are normal the only interesting sheaf is $\alpha_{F_{n}}^{2}$ because $\beta^{1} = \alpha^{1} = \Omega^{1}$ and $\alpha^{2} = \beta^{2}$.

\begin{lemma}\label{$F_{n}$}
For $n = 2p \geq 4$ the form $(a.b)^{p}.da\wedge db/z^{2p-1}$ is a section of $\alpha_{F_{2p}}^{2}$ and the form $(a.b)^{p-1}.da\wedge db/z^{p-1}$ is also a section of $\alpha_{F_{2p}}^{2}$.\\
For $n = 2p +1 \geq 3$ the form $(a.b)^{p}.da\wedge db/z^{2p}$ is a section of $\alpha_{F_{2p+1}}^{2}$.\\
Moreover, all  these forms are not holomorphic forms.
\end{lemma}

\parag{Proof} On $F_{n}$ we have the equalities
\begin{align*}
& a^{n-1}.da\wedge db = z^{n-1}.dz\wedge db \quad {\rm and} \\
& b^{n-1}.da\wedge db = z^{n-1}.dz\wedge da
\end{align*}
so we have
$$ (dz\wedge da).(dz\wedge db) = \frac{(a.b)^{n-1}.(da\wedge db)^{2}}{z^{2n-2}} $$
which implies
$$ \left(\frac{(a.b)^{p}.da\wedge db}{z^{2p-1}}\right)^{2} = a.b.(dz\wedge da).(dz\wedge db) \quad {\rm for} \  n = 2p $$
and
$$  \left(\frac{(a.b)^{p}.da\wedge db}{z^{2p}}\right)^{2} =  a.b.(dz\wedge da).(dz\wedge db) \quad {\rm for} \  n = 2p+1. $$
To prove that $(a.b)^{p-1}.da\wedge db/z^{p-1}$ is  a section of $\alpha_{F_{2p}}^{2}$ consider the map $f : F_{2p} \to S_{2p}$ given by $f(a, b, z) = (a^{p}- b^{p}, a^{p} + b^{p}, z)$ and compute the pull back of the form $dx\wedge dy/z^{p-1}$ which is a section of $\alpha^2_{S_{2p}}$ (see above). The result follows.\\
To see that these forms are not holomorphic is a simple exercise using the homogeneity on $F_{n}$; we live it to the reader.$\hfill \blacksquare$\\

\end{document}